\documentclass[11pt, reqno]{amsart}
\usepackage{lipsum}
\usepackage{a4wide}
\usepackage{amssymb} 
\usepackage{amsmath}
\usepackage{amsthm} 
\usepackage{tikz} 
\usepackage{amscd}
\usepackage{hyperref}
\usepackage{enumitem}
\hypersetup{colorlinks,linkcolor={blue},citecolor={blue},urlcolor={red}}
\theoremstyle{plain}
\usepackage[margin=2.5cm]{geometry}
\numberwithin{equation}{section}

\newcommand{\diag}{\operatorname{diag}}

\newcommand{\End}{\operatorname{End}}

\newtheorem{theorem}{Theorem}[section]
\newtheorem{corollary}[theorem]{Corollary} 

\newtheorem{remark}[theorem]{Remark}

\title{Quotients of commuting schemes associated to Symmetric Pairs}
\author{ Santosh Nadimpalli, Santosha Pattanayak}
\date{\today}
\begin{document}
\maketitle
\begin{abstract}
Let $\mathfrak{g}$ be a classical Lie algebra over an algebraically
closed field $k$ of characteristic zero. Let $\theta$ be an involution
of $\mathfrak{g}$, and let $\mathfrak{g}_0$ and $\mathfrak{g}_1$ be
$1$ and $-1$ eigenspaces of $\theta$. Let
  $G$ be a classical group with Lie algebra $\mathfrak{g}$ and let
  $G_0$ be the connected subgroup of $G$ with
  ${\rm Lie} (G_0)=\mathfrak g_0$. For $d \geq 2$, let
  $\mathfrak{C}^d(\mathfrak{g}_1)$ be the $d$-th commuting scheme
  associated with the symmetric pair $(\mathfrak g, \mathfrak
  g_0)$. In this article, we study the categorical quotient
  $\mathfrak{C}^d(\mathfrak{g}_1)//{G_0}$ via the Chevalley
  restriction map. As a consequence we study the normality and
  reducedness of the scheme
  $\mathfrak C^d(\mathfrak g_1)//G_0$.  As a part
  of the proof, we describe a generating set for the algebra
  $k[\mathfrak{g}_1^d]^{G_0}$, which is of independent interest.
\end{abstract}
\section{Introduction}
Let $\mathfrak{g}$ be a reductive Lie algebra over an algebraically
closed field $k$ of characteristic zero, and let $\tau$ be an
involution of $\mathfrak{g}$. The eigenspaces of $\tau$ determine a
$\mathbb Z_2$-grading of $\mathfrak{g}$, we have
$\mathfrak g=\mathfrak g_0 \oplus \mathfrak g_1$. We say that
$(\mathfrak g, \mathfrak g_0)$ is a symmetric pair (see \cite[Chapter 37]{tauvvel}
and \cite[Chapter X]{helg}
for a general treatment of symmetric pairs). Let $G$ be a reductive
group with Lie algebra $\mathfrak g$ and $G_0$, the connected subgroup
of $G$ with Lie algebra $\mathfrak g_0$. Then $\mathfrak g_1$ is a
$\mathfrak g_0$-module.  Let $\mathfrak c \subset \mathfrak g_1$ be a
maximal subspace consisting of pairwise commuting semisimple
elements. Any such subspace is called a Cartan subspace. All Cartan
subspaces are $G_0$-conjugate.  The number $\dim_k \mathfrak c$ is
called the rank of the symmetric pair. Let $N(\mathfrak c)$ and
$Z(\mathfrak c)$ be the normaliser and the centraliser of
$\mathfrak c$ in $G_0$ respectively. Then $Z(\mathfrak c)$ is a normal
subgroup of $N(\mathfrak c)$ and the quotient
$W_{\mathfrak c} :=N(\mathfrak c)/Z(\mathfrak c)$ is a finite group,
called the little Weyl group of the pair $(G,G_0)$. From \cite{Vin},
it is known that the action of $W_{\mathfrak c}$ on $\mathfrak c$ is
generated by transformations fixing a hyperplane in $\mathfrak c$,
that the restriction map
$k[\mathfrak g_1]^{G_0} \rightarrow k[\mathfrak c]^{W_{\mathfrak c}}$
is an isomorphism.

For $d \geq 2$, the $d$-th commuting scheme
associated with the involution $\tau$, or with the symmetric pair
$(\mathfrak g, \mathfrak g_0)$, is defined
as
$$\mathfrak{C}^d(\mathfrak{g}_1) :=\{(x_1,x_2, \cdots, x_d) \in
\mathfrak{g}_1^d : [x_i,x_j]=0, \ \text{for all}\ i, j\in \{1, 2\dots,
n\}\}.$$ The scheme $\mathfrak{C}^d(\mathfrak{g}_1) $ naturally
generalises the commuting scheme defined by
$$\mathfrak{C}^{d}(\mathfrak{g})=
\{(x_1,x_2, \cdots, x_d) \in
\mathfrak{g}^d : [x_i,x_j]=0, \ \text{for all}\ i, j\in \{1, 2\dots, n\}\}.$$
Richardson in the article \cite{richard}, proved that the scheme
$\mathfrak{C}^2(\mathfrak{g})$ is irreducible.  For $d\geq 4$, the
scheme $\mathfrak{C}^d(\mathfrak{gl}_n)$ is irreducible if and only if
$n\leq 3$ (see \cite{dom_com_mat}). So, the scheme
$\mathfrak{C}^d(\mathfrak{g})$ is not necessarily irreducible, for
large enough $d$.  It is a well known conjecture that the commuting
scheme $\mathfrak C^2(\mathfrak g)$ is reduced. For $d \geq 3$ it is
not expected that $\mathfrak C^d(\mathfrak g)$ is reduced. 

However, the categorical quotient
$\mathfrak C^d(\mathfrak g)//G:={\rm Spec}(k[\mathfrak C^d(\mathfrak
g)]^G)$ behaves better. A weak version of the Chevalley restriction
theorem for the commuting scheme is proved by Hunziker in
\cite{hun}. Hunziker proved that the morphism
$i: \mathfrak h^d//W \rightarrow \mathfrak C^d(\mathfrak g)//G$ of
quotient schemes induced by the inclusion map
$\mathfrak h^d \rightarrow \mathfrak C^d(\mathfrak g)$ is a finite
bijective morphism, where $\mathfrak h$ is a Cartan subalgebra of
$\mathfrak g$. In other words, $i$ is a universal homeomorphism of
schemes.  In particular, $\mathfrak h^d//W$ is the normalisation of
the underlying reduced subscheme
$(\mathfrak C^d(\mathfrak g)//G)^{red}$.  This also follows from a
result of Luna (see \cite{luna}). Since $\mathfrak h^d//W$ is
irreducible, the categorical quotient $\mathfrak C^d(\mathfrak g)//G$
is also irreducible. In \cite{GG}, Gan and Ginzburg proved that the
quotient scheme $\mathfrak C^2(\mathfrak g)//G$ is reduced in the case
where $G={\rm GL}_n(k)$. In \cite{CC1}, Chen and Ng$\rm{\hat{o}}$
conjectured that the morphism
$i: \mathfrak h^d//W \rightarrow \mathfrak C^d(\mathfrak g)//G$ is an
isomorphism of schemes. Since $\mathfrak h^d//W$ is reduced and
normal, the conjecture implies that the categorical quotient
$\mathfrak C^d(\mathfrak g)//G$ is reduced and normal. They proved a
weaker version of the conjecture for $G={\rm GL}_n(k)$ and that was
sufficient to show that the quotient $\mathfrak C^d(\mathfrak g)//G$
is reduced. For $G={\rm GL}_n$, Vaccarino (\cite{vac}) also proved
that the map
$i: \mathfrak h^d//W \rightarrow \mathfrak C^d(\mathfrak g)//G$ is an
isomorphism. Both the proofs use the spectral data map constructed by
Deligne which we recall in section 3. In a recent preprint, Losev
proved that the almost commuting scheme, denoted by $X_n$, associated
to the symplectic Lie algebra is irreducible, reduced complete
intersection; and the categorical quotient
$X_n//{\rm Sp}_{2n}(\mathbb{C})$ is isomorphic to
$\mathfrak{C}^2(\mathfrak{sp}_{2n})/{\rm Sp}_{2n}(\mathbb{C})$ (see
\cite{losev2021commuting}). 

The irreducibility problem for the commuting varieties associated to
symmetric pairs was first considered by Panyushev in \cite{pan1}, and
he observed that $\mathfrak C^2(\mathfrak g_1)$ can be reducible. He
showed that $\mathfrak C^2(\mathfrak g_1)$ is irreducible if and only
if
$\mathfrak C^2(\mathfrak g_1)=\overline{G_0(\mathfrak c \times
  \mathfrak c)}$. Using this he concluded that if
$(\mathfrak g, \mathfrak g_0)$ is a symmetric pair of maximal rank
(i.e., rk($(\mathfrak g, \mathfrak g_0))=rk(\mathfrak g)$), then
$\mathfrak C^2(\mathfrak g_1)$ is an irreducible normal complete
intersection and the ideal of $\mathfrak C^2(\mathfrak g_1)$ in
$k[\mathfrak g_1 \times \mathfrak g_1]$ is generated by quadrics. In
\cite{pan1} Panyushev showed that for a symmetric pair
$(\mathfrak g, \mathfrak g_0)$ of maximal rank, the quotient variety
$\mathfrak{C}^2(\mathfrak{g}_1)//G_0$ is isomorphic to
$\mathfrak h^2//W$, where $\mathfrak h$ is a Cartan subalgebra of
$\mathfrak g$ and $W$ is the Weyl group of $\mathfrak g$. Later (see
\cite{pan2} and \cite{pan3}) he extended the irreducibility results
for some more symmetric pairs including the pairs that we are
interested in this paper. See also \cite{SR1} and \cite{SR2} for some
more results in this direction.

In this note, following the work of Chen and Ng$\rm{\hat{o}}$ in
\cite{CC2} we consider all classical symmetric pairs except the pairs
$(\mathfrak{so}_{2n}(k), \mathfrak{gl}_n(k))$,
$(\mathfrak{sp}_{2(n+m)}(k), \mathfrak{sp}_{2n}(k)\times
\mathfrak{sp}_{2m}(k))$ and
$(\mathfrak{so}_n(k) \times \mathfrak{so}_n(k), \mathfrak{so}_n(k))$
and show that the induced morphism
$\mathfrak i: \mathfrak c^d//W_{\mathfrak c} \rightarrow \mathfrak
C^d(\mathfrak g_1)//G_0$ is an isomorphism of affine schemes. Since
$\mathfrak c^d//W_{\mathfrak c}$ is normal and reduced, the
isomorphism implies that the quotient scheme
$\mathfrak C^d(\mathfrak g_1)//G_0$ is normal and reduced. The main
idea of the proofs is to construct a section (called the spectral data
map) of the Chevalley restriction map using either the determinant map
or the Pfaffian norm map (see Chen and Ng$\rm{\hat{o}}$
\cite{CC2}). In order to verify that the spectral data map is 
a section of Chevalley restriction, we need a convenient set of
generators for the algebra $k[\mathfrak{g}_1^d]^{G_0}$. We use
fundamental theorems of classical
invariant theory to produce a set of generators for
$k[\mathfrak{g}_1^d]^{G_0}$. Then, the explicit form of the restriction
to Cartan subspace 
of these generators already show the surjectivity of Chevalley
restriction map. The injectivity is then verified by formal properties
of the spectral data map. 
\subsection{Classical Symmetric Pairs}
Let $\mathfrak{g}$ be a classical Lie algebra and let $\tau$ be an
involution of $\mathfrak{g}$. Let $\mathfrak{g}_0$ and
$\mathfrak{g}_1$ be the $1$ and $-1$ eigenspaces of $\tau$. The pair
$(\mathfrak{g}, \mathfrak{g}_0)$ is called a classical symmetric pair.
The following is the classification of classical symmetric pairs:
(see for example \cite[Chapter X]{helg})
\begin{enumerate}
\item Bilinear Form:
$(AI) : (\mathfrak{gl}_n(k), \mathfrak{so}_n(k)), \,\, \,\, (AII) :
(\mathfrak{gl}_{2n}(k), \mathfrak{sp}_{2n}(k))$,

\item Polarization: $(DIII): (\mathfrak{so}_{2n}(k), \mathfrak{gl}_n(k))$,
$(CI) : (\mathfrak{sp}_{2n}(k), \mathfrak{gl}_n(k))$,

\item Direct Sum:
  $(AIII): (\mathfrak{gl}_{m+n}(k), \mathfrak{gl}_m(k) \times
  \mathfrak{gl}_n(k))$,
  $(BDI) : (\mathfrak{so}_{m+n}(k), \mathfrak{so}_m(k) \times
  \mathfrak{so}_n(k))$,
  $(CII) : (\mathfrak{sp}_{2(m+n)}(k), \mathfrak{sp}_{2m}(k) \times
  \mathfrak{sp}_{2n}(k))$.
\end{enumerate}
See Section 5 for the detailed structure of these symmetric pairs. 
The main theorem of this paper is the following:
  \begin{theorem}\label{main}
   Let $(\mathfrak{g}, \mathfrak{g}_0)$ be a classical symmetric pair
except the pair
$(\mathfrak{so}_{2n}(k), \mathfrak{gl}_n(k))$.  The Chevalley restriction map
    $\mathfrak i: \mathfrak{c}^d//W_\mathfrak{c} \rightarrow
    \mathfrak{C}^d(\mathfrak{g}_1)//{G_0}$ is an isomorphism of affine
    schemes.
  \end{theorem}
  Our techniques for the construction of a section for the Chevalley
  restriction map does not work in the case where
  $(\mathfrak{g}, \mathfrak{g}_0)$ is equal to
  $(\mathfrak{so}_{2n}(k), \mathfrak{gl}_n(k))$. The reason being that
  the rank of symmetric pair is half of the rank of rank of
  $\mathfrak{g}_0$. See remark \ref{except} for details.

  Since
  $\mathfrak c^d//W_{\mathfrak c}$ is normal and reduced, as an
  immediate corollary we get the following:
  \begin{corollary}
  The categorical quotient scheme $\mathfrak C^d(\mathfrak g_1)//G_0$ is
  normal and reduced.
  \end{corollary} 
  
  We will explain the contents of each section. In section 2, we
  recall the notion of polynomial laws and Roby's results. In section
  3, we recall Deligne's construction of the spectral data map. In
  Section 4, we describe the ring of $G_0$ invariants of
  $ k[\mathfrak{g}_1^d]$. For many pairs
  $(\mathfrak{g}, \mathfrak{g}_1)$ these results are new and of
  independent interest. In section 5, we give a proof of the main
  theorem by constructing the respective spectral data map for each
  symmetric pair and verify that this spectral data map is a section
  of the Chevalley restriction map. 
\section{Polynomial Laws}
\subsection{}
Let $A$ be a commutative ring and let ${\rm Alg}_A$ be the category of
$A$-algebras. Let $V$ be an $A$-module and let $V_A$ be the functor
from ${\rm Alg}_A$ to the category of $A$-modules given by
$R\mapsto V\otimes _A R$. For two $A$-modules $V$ and $W$, we denote
by $P(V, W)$, the set of natural transformations between the functors
$V_A$ and $W_A$. If $W=A$, then we denote by $P(V)$ the set $P(V,
A)$. The set $P(V)$ is called the set of {\it polynomial laws} on
$V$. Let $S_A$ be the polynomial ring $A[X_1, X_2,\dots, X_n]$. For
$f\in P(V)$ we get a map $f_{S_A}:V\otimes_AS_A\rightarrow S_A$. Thus
for any $f\in P(V)$, and a finite set of elements $v_1,\dots, v_n$, we
associate a polynomial $P_f\in S_A$ given by
$f_{S_A}(v_1X_1+v_2X_2+\dots+v_nX_n)$. A polynomial law $f\in P(V)$ is
called homogenous of degree $d$ if $f(uv)=u^nf(v)$, for all $A$-algebras $R$, 
for all $u\in R^\times$ and $v\in V\otimes _AR$. We denote by $P_n(V)$ the set
of all degree $n$ homogenous polynomial laws on $V$. For a general
reference for polynomial laws we refer to \cite{Roby}.
\subsection{}
Let $V$ be an $A$-module, and let $T^n(V)$ be the  $n$-fold
tensor product of the $A$-module $V$. We denote by $TS^n(V)$, the
$A$-submodule fixed by the action of $S_n$ on $T^n(V)$ and by
$S^n(V)$, the $S_n$-coinvariants of $T^n(V)$.
Roby (see \cite{Roby}) showed
that the homogenous polynomial laws on $V$ of degree $n$ are in
canonical bijection with homogenous degree $1$ polynomial laws on
$TS^n(V)$ given by the relation
$$f(v)=h(v^{\otimes n}), f\in P_n(V), h\in P_1(TS^n(V)),$$ where $P_1(TS^n(V))$ 
is the space of degree $1$ polynomial laws on $TS^n(V)$.  Moreover, if
$V$ is an $A$-algebra, which is free as an $A$-module, and $f$ is a
multiplicative homogenous polynomial law of degree $n$ on $V$, then
the degree $1$ polynomial law $h$ associated to $f$ is a homomorphism of
algebras
$$TS^n(V)\rightarrow A.$$
\section{Spectral data map}
In this section, we first recall Deligne's construction of spectral
data map for ${\rm GL}_n(k)$ (see \cite[Section
6.3.1]{SGA4_cohomo_cpt_supp}) which plays a pivotal role in the proof
of the main theorem.  As the name suggests it assigns the set of
common eigenvalues of a set of commuting matrices. In fact, we
construct the spectral data map for the relevant symmetric pairs and
show that it is actually a section of the Chevalley restriction
map. For convenience, we recall the Chevalley restriction map
here. Let $G$ be a reductive group over $k$ with Lie algebra
$\mathfrak g$. Let $\mathfrak h$ be a Cartan subalgebra of
$\mathfrak g$ and let $W$ be the Weyl group of $\mathfrak g$. The
diagonal adjoint action of $G$ on $\mathfrak g^d$ leaves the commuting
scheme $\mathfrak C^d(\mathfrak g)$ invariant. The inclusion
$\mathfrak h^d \rightarrow \mathfrak g^d$ factors through the
commuting scheme $\mathfrak C^d(\mathfrak g)$ and it induces a
homomorphism of $k$-algebras
$\mathfrak i: k[\mathfrak C^d(\mathfrak g)]^G \rightarrow k[\mathfrak
h^d]^W$ because the restriction of a $G$-invariant function to
$\mathfrak h^d$ is also $W$-invariant. In other words, we have a
morphism of affine schemes (which we still denote by $\mathfrak i$):
$$\mathfrak i: \mathfrak h^d//W \rightarrow \mathfrak C^d(\mathfrak
g)//G.$$

\subsection{} Let $V$ be an $n$-dimensional $k$-vector space and and
let $R$ be any $k$-algebra. Let $\mathfrak{C}^d(\mathfrak{gl}(V))$ be
the pair wise commuting scheme of $\mathfrak{gl}(V)$ whose $R$-points
are given by the set
$$\{(x_1, \dots, x_d)\in \mathfrak{gl}(V\otimes R)^d : [x_i,
x_j]=0,\ \text{for all}\ i, j\in [d]\}.$$ Let $A$ be a $k$-algebra
representing the functor $\mathfrak{C}^d(\mathfrak {gl}(V))$. Let
$(x_1, x_2,\dots, x_d)$ be the universal point in
$\mathfrak{C}^d(A)$. We then have a homomorphism of $R$-algebras:
$$p:k[X_1, \dots, X_d]\otimes R\rightarrow \mathfrak{gl}(A\otimes R)
$$
given by $p(X_1,\dots, X_d)=(x_1, \dots, x_d)$. The map $\det\circ p$
gives a homogenous polynomial law on the algebra $k[X_1, \dots,
X_d]$. Since $\det\circ p$ is multiplicative, using Roby's result we
get an algebra homomorphism 
$$\tilde{\mathfrak s}:(k[X_1, X_2,\dots X_d]^{\otimes n})^{S_n}\rightarrow A$$
such that $\det p((f))=\tilde{\mathfrak s}(f^{\otimes n})$ for all
$f \in k[X_1, X_2,\dots X_d]$. We note that the image of the map
$\tilde{\mathfrak s}$ belongs to the algebra of $GL(V)$-invariants
$A^{{\rm GL}(V)}$. Thus we obtain a map of schemes
$$\mathfrak s:
\mathfrak{C}^d(\mathfrak{gl}(V))//{\rm GL}(V) \rightarrow
\mathfrak{h}^d//{S_n}.$$ Here, $\mathfrak{h}$ is a Cartan subalgebra
of $\mathfrak{gl}(V)$. This map $\mathfrak s$ is called the spectral data
map.  Note that in this case the Weyl group is $S_n$ and the Cartan
subalgebra is the standard representation of $S_n$. The compositions
$\mathfrak s\circ \mathfrak i$ and $\mathfrak i\circ \mathfrak s$ are
verified to be identity on a generating set given by trace functions
due to Procesi (see \cite{procesi}).  As a result it was concluded
that the quotient scheme $\mathfrak{C}^d(\mathfrak{gl}(V))//GL(V)$ is
normal and reduced (see \cite{CC1} for details).
\section{Invariants}
Let $(\mathfrak{g}, \mathfrak{g}_0)$ be a classical symmetric pair so
that $\mathfrak g=\mathfrak g_0 \oplus \mathfrak g_1$. We require a
generating set for the $k$-algebra $k[\mathfrak{g}_1^d]^{G_0}$. In
some of the cases, these invariants are known from the work of
Procesi (see \cite{procesi}). In this section, we consider those cases which are not
available (to the best of our knowledge) in the literature.  These
invariants will be used to check that the spectral data map is a
section of the Chevalley restriction map.
\subsection{}\label{invariants_gl_n_gl_m}
Let $W(n, m, d)$ be the space $M_{n\times
   m}^d(k)\oplus M^d_{m\times n}(k)$. Let $G$ be the
 group ${\rm GL}_n(k)\times {\rm GL}_m(k)$ via left
 and right action on $W(n, m,d)$. We define
 $$\mu: W(n, m, d)\rightarrow
 M^{d^2}_{n\times n}(k); ((M_1, \dots, M_d), (N_1,\dots, N_d)) \mapsto
 (A_{ij}),$$ where $A_{ij}=M_iN_j$. Note that the map
 $\mu^\ast: k[M^{d^2}_{n\times n}(k)]\rightarrow k[M_{n\times
   m}^d(k)\bigoplus M^d_{m\times n}(k)]^{{\rm GL}_m(k)}$ is surjective
 from the first fundamental theorem for $GL_n(k)$ (see Theorem 5.2.1 of
 \cite{GW}). Hence, we get that the map
 \begin{equation}\label{w(n,m,d)}
   \mu^\ast:
 k[M^{d^2}_{n\times n}(k)]^{{\rm GL}_n(k)}
 \rightarrow
 k[M_{n\times m}^d(k)\oplus
 M^d_{m\times n}(k)]^G
\end{equation} is surjective. Since $k[M^{d^2}_{n\times n}(k)]^{{\rm GL}_n(k)}$ 
is generated by monomials in $A_i$, $i \in [d^2]$ (see Theorem 3.4 of \cite{procesi}), 
from \eqref{w(n,m,d)} we get that polynomials of the form
 $${\rm Tr}(M_1M_2\dots M_l)$$
 where $M_{i}=Q_{n_i}R_{m_i}$ with $Q_{n_i}, R_{m_i}^t\in M_{n\times
 m}(k)$ and $n_i, m_i\in [d]$ for all $i\in [l]$ generate the algebra 
$k[W(n,m,d)]^G$. 
\subsection{}\label{invariants_so_so}
Let $V(n, m, d)$ be the space $M^d_{n\times m}(k)$,
and let $G$ be the group ${\rm SO}_n(k)\times
{\rm SO}_m(k)$. The group $G$ acts on $V(m,n, d)$ by setting 
$$(A, B)(M_1, M_2,\dots, M_d)=(AM_1B^t, AM_2B^t,\dots, AM_dB^t).$$
Let $\mu:M^d_{n\times m}(k)\rightarrow
({\rm Sym}^2k^m)^{d}\bigoplus_{i>j}M_{m\times m}(k)$ be the map
$$(M_1, M_2,\dots, M_d)\mapsto (A_{ij}: i\geq j),$$
where $A_{ij}=M_i^tM_j$.  The map
$\mu^\ast:k[({\rm Sym}^2k^m)^{d}\bigoplus_{i>j}M_{m\times
  m}(k)]\rightarrow k[M^d_{n\times m}(k)]^{{\rm O}_n(k)}$ is a
surjective map from the first fundamental theorem for orthogonal group
(see Theorem 5.2.2 of \cite{GW}). Since ${\rm SO}_m(k)$ is reductive, we get
that the map
\begin{equation}\label{surj}
  \mu^\ast:k[({\rm Sym}^2k^m)^{d}\bigoplus_{i>j}M_{m\times
  m}(k)]^{{\rm
      SO}_m(k)}
  \rightarrow
  k[M^d_{n\times m}(k)]^{{\rm O}_n(k)\times
SO_m(k)}
\end{equation} is surjective. 
The algebra $k[({\rm Sym}^2k^m)^{d}\bigoplus_{i>j}M_{m\times
  m}(k)]^{{\rm
    {\rm SO}}_m(k)}$ is spanned by polynomials of the form
$${\rm Tr}(A_{n_1m_1}A_{n_2m_2}\dots
A_{n_lm_l}),$$
and
$$\widetilde{{\rm Pf}}(A_{n_1m_1}A_{n_2m_2}\dots
A_{n_{l'}m_{l'}}),$$ where $\widetilde{Pf}$ is the complete polarisation of the
pfaffian, and $A_{ii}$ is a symmetric matrix and $A_{ij}$ for $i\neq j$
is any $m\times m$ matrix; for details see \cite{ATZ}. Hence, the algebra
$$k[M^d_{n\times m}(k)]^{{\rm O}_n(k)\times
SO_m(k)}$$
is generated by polynomials of the form 
\begin{align}\label{odd_even_shit}
 &{\rm Tr}(A_{n_1m_1}A_{n_2m_2}\dots
A_{n_lm_l}),\\
 \label{paffian}&  \widetilde{{\rm Pf}}(A_{n_1m_1}A_{n_2m_2}\dots
A_{n_l'm_l'}),
\end{align}
where $A_{n_im_i}=M_{n_i}^tM_{m_j}$. If $n$ is odd then we get that
a system of generators of
$k[M^d_{n\times m}(k)]^{{\rm SO}_n(k)\times {\rm SO}_m(k)}$ is given
by the invariants in \eqref{odd_even_shit}.
\begin{remark}\normalfont
  If $M_{i}^tM_j=M_j^tM_i$ and $M_jM_i^t=M_iM_j^t$ for all
  $i, j\in [d]$. Then the matrix
  $$M_{n_1}^tM_{m_1}M_{n_2}^tM_{m_2}\dots M_{n_l}^tM_{m_l}$$ is
  symmetric. Hence the restriction of the invariant \eqref{paffian} to
  the subvariety of $M^d_{n\times m}(k)$ defined by the relations
  $M_{i}^tM_j=M_j^tM_i$ and $M_jM_i^t=M_iM_j^t$ is zero. 
  \end{remark}

If $n$ is even and $n=m$,
then the algebra
$k[M^d_{n\times m}(k)]^{{\rm SO}_n(k)\times {\rm SO}_m(k)}$ might be
strictly bigger than
$k[M^d_{n\times m}(k)]^{{\rm O}_n(k)\times {\rm SO}_m(k)}$; for
instance the elements
$$\det(T_1\otimes A_1+T_2\otimes A_2+\cdots+T_d\otimes A_d)$$
for any $(T_1, T_2,\dots, T_d)\in M^d_{l\times l}(k)$, where $l \geq 1$, are invariant
for ${\rm SO}_n(k)\times {\rm SO}_m(k)$.

  
\subsection{}\label{invariants_sym_gl_n}
Let $W_{(d_1, d_2)}$ be the space
$({\rm Sym}^2(V))^{d_1}\oplus ({\rm Sym}^2(V^\ast))^{d_2}$, and let
$G$ be the group ${\rm GL}(V)$. For the definitions of full
polarization and restitution, we refer to Section 3.2.2 of
\cite{procesi-book}. The algebra $A=k[W_{(d_1, d_2)}]^{G}$ is graded
by $\mathbb{N}^{d_1+d_2}$. Let $F$ be a non-zero multihomogeneous
invariant in $A$. Let $P(F)$ be the full polarization of $F$, and note
that $P(F)$ is a multilinear invariant on $W_{d'_1, d_2'}$ for some
integers $d'_1$ and $d_2'$ depending on the degree of $F$. We embed
$W_{d'_1, d_2'}$ in
$W_{d_1', d_2'}':=(V\otimes V)^{d_1'}\oplus (V^\ast\otimes
V^\ast)^{d_2'}$. The space of multilinear invariants of
$W_{d_1', d_2'}'$ is equal to
$(V^{\otimes 2d_1'}\otimes V^{\otimes 2d_2'})^{G}$. We then have
$d_1'=d_2'=d'$ and the space $(V^{2d'}\otimes V^{2d'})^{G}$ is spanned
by the complete contractions (see Corollary 5.3.2 of \cite{GW}). So
$(V^{2d'}\otimes V^{2d'})^{G}$ is spanned by monomials of the form
$$u_\sigma=\prod_{i=1}^{2d'}u_{i, \sigma(i)},$$
where $\sigma\in S_{2d'}$ and 
$$u_{ij}(v_1\otimes v_2\otimes\cdots\otimes v_{2d'}\otimes
v_1^\ast\otimes v_2^\ast\otimes \cdots\otimes
v_{2d'}^\ast)=v_j^\ast(v_i).$$ Since $G$ is reductive, the restrictions 
of $u_\sigma$ from $W'_{d', d'}$ to the space of multilinear
invariants of $W_{d', d'}$ generate the multilinear $G$-invariants. Since every 
multihomogeneous invariant is the restitution of some multilinear invariant, the 
restitution of $u_\sigma$, $\sigma \in S_{2d'}$ generate $k[W_{(d', d')}]^{G}$. 
The restitution of $u_\sigma$ is of the form
$${\rm Tr}(M_1M_2\dots M_{k}),$$
with $M_i$ is of the from $Q_{n_i}R_{m_i}$, where $(Q_1,\dots, Q_d)$
and $(R_1,\dots, R_d)$ are tuples of symmetric matrices and
$n_i, m_i\in [d]$ (see \cite[Exercise 1 of 5.3.3]{GW}).

  \subsection{}\label{invariants-type-cii}
  Let $V(2p, 2q,d)$ be the space $M^d_{2p\times 2q}(k)$ and let $G$ be
  the group ${\rm Sp}_{2p}(k)\times {\rm Sp}_{2q}(k)$. We define the
  action of $G$ on $V(2p, 2q,d)$  by setting:
  $$(A, B)(M_1, M_2,\dots, M_d)=(AM_1B^{-1}, AM_2B^{-1},\dots,
  AM_dB^{-1}).$$ Let $T_r$ be the $2r\times 2r$ matrix given by
$\diag(\mu, \mu,\dots, \mu)$, where $\mu$ is the
matrix $$\begin{pmatrix}0&1\\-1&0\end{pmatrix}.$$
  Let $\mu:V(2p, 2q, d)\twoheadrightarrow M^{d^2}_{2q, 2q}(k)$ be given by
  $(M_1, M_2\dots, M_d)\mapsto (A_{ij})$, where
  $A_{ij}=M_i^tT_{p}M_jT_q$. The map $\mu$ induces a surjective map
 $\mu^\ast:k[M^{d^2}_{2q\times 2q}(k)]\twoheadrightarrow k[V(2p, 2q,
 d)]^{{\rm Sp}_{2p}(k)}$. Hence, we get a surjective map
 $$\mu^\ast:k[M^{d^2}_{2q\times 2q}]^{{\rm Sp}_{2q}(k)}\twoheadrightarrow
 k[V(2p, 2q, d)]^{G}.$$
 Note that the algebra $k[M^{d^2}_{2q\times 2q}(k)]^{{\rm Sp}_{2q}(k)}$ is
 spanned by polynomials:
 $${\rm Tr}(M)$$
 where $M$ is a monomial in $A_{ij}$ or $A_{ij}^t$ with $A_{ij}$ is a
 variable matrix of size $2q\times 2q$, for all $i,j\in[d]$ (see
 Theorem 10.1 of \cite{procesi}). Hence, the $k$-algebra
 $k[V(2p, 2q, d)]^{G}$ is spanned by elements of the form
 $${\rm Tr}(M_{n_1m_1}M_{n_2m_2}\dots M_{n_lm_l})$$
 where $M_{ij}=M_jT_{p}M_i^tT_q$ and $n_i, m_i\in [d]$. 
\section{Construction of Spectral data map and proof of the main
  theorem}
In this section we construct the spectral data map for the relevant
classical symmetric pairs and we show that in each case the spectral
data map is a section of the Chevalley restriction map.
\subsection{The symmetric pair \texorpdfstring{$AI$}{}}

Let $V$ be a $k$-vector space of dimension $n$, and let
$\omega:V\times V\rightarrow k$ be a non-degenerate symmetric bilinear
form on $V$.  Let
$\tau: \mathfrak{gl}(V) \rightarrow \mathfrak{gl}(V)$ be the
involution defined by $\tau(X)=-X^\ast$, where $X^\ast$ is given by:
$$\omega(Xv, w)=\omega(v, X^\ast w).$$
Let $\mathfrak{g}_0$ and $\mathfrak{g}_1$ be the $1$ and $-1$
eigenspaces of $\tau$. Then the algebra $\mathfrak{g}_0$ can be
identified with the Lie algebra of the connected component $G_0$ of
${\rm O}(V, \omega)$ containing the identity element.  With respect to
an orthogonal basis $(v_1, v_2,\dots, v_n)$ of $V$ a Cartan subspace
$\mathfrak{c} \subset \mathfrak{g}_1$ is given by the following
subspace:
$$\{\diag(b_{1}, b_{2},\dots, b_{n}): b_i\in k, i\in [n]\}.$$
The little Weyl group $W_\mathfrak{c}$ is isomorphic to $S_n$--which is
identified with the group of permutations of the
coordinates of $\mathfrak{c}$. 

Let $\mathfrak{C}^d(\mathfrak{g}_1)$ be the
$d$-fold commuting scheme associated with the symmetric pair
$(\mathfrak{gl}(V), \mathfrak{g}_0)$.  We have the Chevalley restriction map
$$\mathfrak i: \mathfrak{c}^d//W_\mathfrak{c}
\rightarrow \mathfrak{C}^d(\mathfrak{g}_1)//{G_0},$$ which is induced from
the inclusion $\mathfrak c^d \rightarrow \mathfrak g_1^d$.  We will
construct a section of this map. 

Let $R$ be any $k$-algebra, and let $A$ be the coordinate ring of
$\mathfrak{C}_d(\mathfrak{g}_1)$.  Let $(x_1, x_2,\dots, x_d)$ be
the universal point of 
$\mathfrak{C}_d(\mathfrak{g}_1)(A)$. For any $k$-algebra $R$, we
define the map
$$p:R[X_1,X_2,\dots, X_d]\rightarrow
\mathfrak{gl}_n(A\otimes R), \ X_i\mapsto x_i, i\in [d].$$
Note that the map $\det\circ\ p$ is a degree $n$ multiplicative
map. Hence, by Roby's theorem we get an algebra map
$$\tilde{\mathfrak{s}}:TS^n(k[X_1, X_2,\dots, X_d])\rightarrow A$$
such that
$$\tilde{\mathfrak{s}}(\theta^{\otimes n})=\det\circ\ p(\theta)$$ for
all $\theta\in k[X_1, X_2,\dots, X_d]$.  Since $\det$ is
$G_0$-invariant, the image of the map $\tilde{\mathfrak s}$ is is contained in
$A^{G_0}$.

 Let $(y_1,\dots, y_d)$ be the
tautological point of $\mathfrak{c}(B)$, where
$B=k[\mathfrak{c}^d]$. Note that $B=k[\mathfrak{c}^d]$ is a polynomial
algebra in the variables $b_j(y_i)$, $1 \leq i \leq d$ and
$1 \leq j \leq n$. Let
  $$\beta: TS^n(k[X_1, X_2,\dots, X_d])\rightarrow k[\mathfrak{c}^d]$$
  be the isomorphism given by $\beta(X_{i, j})=b_j(y_i)$, where
  $X_{i, j}$ is the variable $X_i$ in the $j$-th copy, for $i\in [d]$
  and $j\in [n]$. We get a map
  \begin{equation}\label{section_A1}
    \mathfrak{s}:
    k[\mathfrak{c}^d]^{W_\mathfrak{c}}\rightarrow A^{G_0}
    \end{equation}
  such that $\mathfrak{s}\circ \beta(\theta^{\otimes n})=\det\circ
  p(\theta)$.  Thus we obtain a map of schemes
$\mathfrak s: \mathfrak{C}^d(\mathfrak{g}_1)//{G_0} \rightarrow
\mathfrak{c}^d//W_\mathfrak{c}$. We show that the map
  $\mathfrak s: \mathfrak{C}^d(\mathfrak{g}_1)//{G_0} \rightarrow
  \mathfrak{c}^d//W_\mathfrak{c}$ is the inverse of the Chevalley restriction
  map
  $\mathfrak i: \mathfrak{c}^d//W_\mathfrak{c} \rightarrow
  \mathfrak{C}^d(\mathfrak{g}_1)//{G_0}$.

  For the diagonal action of $G_0$ on the $d$-copies of
  $\mathfrak{gl}(V)$, the ring of invariants
  $k[\mathfrak{gl}(V)^d]^{G_0}$ is generated by the elements $Tr(M)$
  and some polarized Pfaffians, where $M$ is a monomial in
  $X_j, X_j^*$, $j=1,2, \cdots, d$ (see \cite{ATZ}). Recall that the
  space $\mathfrak{g}_1$ is the $-1$ eigen space of the map
  $\theta(X)=-X^*$. So the elements in $\mathfrak{g}_1$ satisfy
  $X^*=X$ and hence the polarized Pfaffians become zero when we
  restrict them to $\mathfrak{g}_1^d$.
 
  Since $\mathfrak{C}^d(\mathfrak{g}_1)$ is a closed subscheme of
  $\mathfrak{gl}(V)^d$ and $G_0$ is (linearly) reductive, the
  restriction map
  $k[\mathfrak{gl}(V)^d]\rightarrow k[\mathfrak{C}^d(\mathfrak{g}_1)]$
  induces a surjective map
  $k[\mathfrak{gl}(V)^d]^{G_0}\rightarrow
  k[\mathfrak{C}^d(\mathfrak{g}_1)]^{G_0}$. So the restrictions of the
  functions $Tr(M)$, where $M$ is a monomial in $X_j$,
  $j=1,2, \cdots, d$ form a generating set for
  $k[\mathfrak{C}^d(\mathfrak{g}_1)]^{G_0}$. Let $(x_1,x_2, \cdots, x_d)$
  be the universal point of the commuting scheme
  $\mathfrak{C}^d(\mathfrak{g}_1)(A)$. Then the set
  $\{ \phi_{\underline{a}}={\rm Tr}(x_1^{a_1}x_2^{a_2}\dots
  x_d^{a_d}): \underline{a}=(a_1,\dots, a_d)\in \mathbb{Z}^d_{\geq
    0}\}$ is a generating set for
  $k[\mathfrak{C}^d(\mathfrak{g}_1)]^G$.
  
 Note that
 $$\mathfrak i(\phi_{\underline{a}})=\sum_{j=1}^n\prod_{i}b_j(y_i)^{a_i}.$$
 From the above equation, we get that
 $\mathfrak{i}(\phi_{\underline{a}})$ generate the $k$-algebra
 $k[\mathfrak{c}^d]^{W_\mathfrak{c}}$. We set
 $\psi_{\underline{a}}=\mathfrak{i}(\phi_{\underline{a}})$.  Let
 $\theta_{\underline{a}} \in k[t]\otimes[X_1,X_2, \cdots ,
 X_d]$ be the polynomial $t-X_1^{a_1}X_2^{a_2}\dots
 X_n^{a_n}$. Then $p(\theta_{\underline{a}})= tI-x_1^{a_1}x_2^{a_2}\dots x_d^{a_d}
 \in \mathfrak g_1(A \otimes k[t])$, where $I$ is the identity matrix.
  Using \eqref{section_A1} we have
  $$\det(p(\theta_{\underline{a}}))=
  \mathfrak s(\beta(\theta_{\underline{a}}^{\otimes n})).$$
By the definition of the map $\beta$ we have
$\beta(\theta_{\underline{a}}^{\otimes
  n})=\prod_{j=1}^n(t-\prod_{i=1}^db_j(y_i)^{a_i})$.
Hence, we get that
\begin{equation}\label{bla_A1}
  \det(tI-x_1^{a_1}x_2^{a_2}\dots x_d^{a_d})
  =\mathfrak s(\prod_{j=1}^n(t-\prod_{i=1}^db_j(y_i)^{a_i})).
 \end{equation}  
  Note that $\det(tI-x_1^{a_1}x_2^{a_2}\dots x_d^{a_d})$ is the
 characteristic polynomial of the matrix
 $x_1^{a_1}x_2^{a_2}\dots x_d^{a_d}$ and hence
 \begin{equation}\label{bla_A2}
   \mathfrak{s}(\det(tI-x_1^{a_1}x_2^{a_2}\dots x_d^{a_d}))=
 t^{n}-\mathfrak{s}({\rm Tr}(x_1^{a_1}x_2^{a_2}\dots
 x_d^{a_d}))t^{n-1}+\ \text{lower degree terms in} \ t.
\end{equation}
Comparing
 the coefficients of $t^{2n-1}$ in equations  \eqref{bla_A1} and
 \eqref{bla_A2} we get
 that $\mathfrak s(\mathfrak
 i(\phi_{\underline{a}}))=\phi_{\underline{a}}$.
Thus, we conclude that $\mathfrak{s}\circ \mathfrak{i}$ and
$\mathfrak{i}\circ \mathfrak{s}$ are both identities. 

 Hence, we get that the Chevalley restriction map
 $\mathfrak{i}:\mathfrak{c}^d//W_\mathfrak{c}\rightarrow
 \mathfrak{C}^d(\mathfrak{g}_1)//G$ is an isomorphism.

\subsection{The symmetric pair \texorpdfstring{$AII$}{}}
Let $V$ be a $k$-vector space of dimension $2n$, and let
$\omega:V\times V\rightarrow k$ be a non-degenerate skew symmetric
bilinear form on $V$.  Let
$\tau: \mathfrak{gl}(V) \rightarrow \mathfrak{gl}(V)$ be the
involution defined by $\tau(X)=-X^\ast$, where $X^\ast$ is given by:
$$\omega(Xv, w)=\omega(v, X^\ast w).$$
Let $\mathfrak{g}_0$ and $\mathfrak{g}_1$ be the $1$ and $-1$
eigenspaces of $\tau$. Then the algebra $\mathfrak{g}_0$ can be
identified with the Lie algebra of the connected component $G_0$ of
${\rm Sp}(V, \omega)$ containing the identity. Let
$$(w_{-n}, w_{-n+1},\dots, w_{-1},w_1, \cdots ,w_{n-1}, w_{n})$$ be a Witt basis of
  $V$, and in this basis a Cartan subspace
  $\mathfrak{c} \subset \mathfrak{g}_1$ is given
  by the following subspace:
$$\{\diag(b_{n},\dots ,b_{1}, b_{1},\dots, b_{n}): b_i\in k, i\in [n]\}.$$
The little Weyl group $W_\mathfrak{c}$ is isomorphic to $S_n$-which is
identified with the group of permutations of the coordinates of
$\mathfrak{c}$.

In \cite{CC2} for any $k$-algebra $R$, the authors
defined a map
\begin{equation}\label{paffian_map_AII}
  N_+: \mathfrak g_1(R) \rightarrow R,
\end{equation}
 called it the
Pfaffian norm map which satisfies that $\det(x)=N_+(x)^2$ for any
$x \in \mathfrak g_1(R)$.  More interestingly, they showed that the
map $N_+$ is multiplicative on the coordinate ring of the commuting
subscheme $\mathfrak C({\mathfrak g_1})$ of
$\mathfrak g_1 \times \mathfrak g_1$.

Let $A$ be the coordinate ring of $\mathfrak{C}^d(\mathfrak{g}_1)$,
and let $(x_1, x_2,\dots, x_d)$ be the universal point of this
scheme. For any $k$-algebra $R$, we define the map
$$p: R[X_1, X_2,\dots, X_d]\rightarrow
\mathfrak{gl}_{2n}(A \otimes R), X_i\mapsto x_i, i\in [d].$$ The image
of $p$ is contained in $\mathfrak{g}_1(A\otimes R)$. Then we have the map
\begin{equation} 
R[X_1, X_2,\dots, X_d]\xrightarrow{p} \mathfrak{g}_1(A\otimes
  R)\xrightarrow{N_+} A \otimes R.
\end{equation}
Since $N_+$ is multiplicative, the composition $N_+ \circ p$ is
multiplicative and homogeneous of degree $n$. Thus by Roby's theorem
we get a homomorphism of $k$-algebras:
$$\tilde{\mathfrak s}:TS^n(k[X_1,\dots,X_d])\rightarrow A$$
such that 
\begin{equation}
  \tilde{\mathfrak s}(q^{\otimes n})=
  N_+(p(\theta)), \,\, \text{for all} \,\, \theta \in k[X_1,\dots, X_d]. 
\end{equation}
Since $N_+$ is a $G_0$-invariant map, the map $\tilde{\mathfrak s}$ is
also $G_0$-invariant and hence we get that the image of
$\tilde{\mathfrak s}$ is contained in $A^{G_0}$.

Let $(y_1,\dots, y_d)$ be the
tautological point of $\mathfrak{c}(B)$, where
$B=k[\mathfrak{c}^d]$. Note that $B=k[\mathfrak{c}^d]$ is a polynomial
algebra in the variables $b_j(y_i)$, $1 \leq i \leq d$ and
$1 \leq j \leq n$. Let
  $$\beta: TS^n(k[X_1, X_2,\dots, X_d])\rightarrow k[\mathfrak{c}^d]$$
  be the isomorphism given by $\beta(X_{i, j})=b_j(y_i)$, where
  $X_{i, j}$ is the variable $X_i$ in the $j$-th copy, for $i\in [d]$
  and $j\in [n]$. Hence, we get a map
  $\mathfrak{s}:k[\mathfrak{c}^d]\rightarrow A^{G_0}$ such that
  \begin{equation}\label{section_AII}
    \mathfrak{s}\circ\beta(\theta^{\otimes n})=N_+\circ p(\theta),
  \end{equation}
  for
  all $\theta\in k[X_1,\dots, X_d]$. Thus we obtain a map of schemes
  $\mathfrak s: \mathfrak{C}^d(\mathfrak{g}_1)//{G_0} \rightarrow
  \mathfrak{c}^d//W_\mathfrak{c}.$ We show that the map
  $\mathfrak s: \mathfrak{C}^d(\mathfrak{g}_1)//{G_0} \rightarrow
  \mathfrak{c}^d//W_\mathfrak{c}$ is the inverse of the Chevalley restriction
  map
  $\mathfrak i: \mathfrak{c}^d//W_\mathfrak{c} \rightarrow
  \mathfrak{C}^d(\mathfrak{g}_1)//{G_0}$.

  For the diagonal action of the symplectic group ${G_0}$ on the
  $d$-copies of $\mathfrak{gl}(V)$, the ring of invariants
  $k[\mathfrak{gl}(V)^d]^{G_0}$ is generated by the elements $Tr(M)$,
  where $M$ is a monomial in $X_j, X_j^*$, $j=1,2, \cdots, d$ of
  degree less than or equal to $2^n-1$ (see Theorem 10.1 of \cite{procesi}). Recall
  that the space $\mathfrak{g}_1$ is the $-1$ eigen space of the map
  $\theta(X)=-X^*$. So the elements in $\mathfrak{g}_1$ satisfy
  $X^*=X$. Since $\mathfrak{C}^d(\mathfrak{g}_1)$ is a closed subscheme of
  $\mathfrak{gl}(V)^d$ and $G_0$ is (linearly) reductive, the
  restriction map
  $k[\mathfrak{gl}(V)^d]\rightarrow k[\mathfrak{C}^d(\mathfrak{g}_1)]$
  induces a surjective map
  $k[\mathfrak{gl}(V)^d]^{G_0}\rightarrow
  k[\mathfrak{C}^d(\mathfrak{g}_1)]^{G_0}$. So the restrictions of the
  functions $Tr(M)$, where $M$ is a monomial in $X_j$,
  $j=1,2, \cdots, d$ form a generating set for
  $k[\mathfrak{C}^d(\mathfrak{g}_1)]^{G_0}$. Let $(x_1,x_2, \cdots, x_d)$
  be the tautological point of the commuting scheme
  $\mathfrak{C}^d(\mathfrak{g}_1)$. Then the set
  $\{ \phi_{\underline{a}}={\rm Tr}(x_1^{a_1}x_2^{a_2}\dots
  x_d^{a_d}): \underline{a}=(a_1,\dots, a_d)\in \mathbb{Z}^d_{\geq
    0}\}$ is a generating set for $k[\mathfrak{C}^d(\mathfrak{g}_1)]^{G_0}$.
  
  Note that
  $$\mathfrak
  i(\phi_{\underline{a}})=2\sum_{j=1}^n\prod_{i}b_j(y_i)^{a_i}.$$
  From the above equation, we get that
 $\mathfrak{i}(\phi_{\underline{a}})$ generate the $k$-algebra
 $k[\mathfrak{c}^d]^{W_\mathfrak{c}}$. 
  We set $\psi_{\underline{a}}=\mathfrak i(\phi_{\underline{a}})$, for
  all $n$-tuples $\underline{a}\in \mathbb{Z}^n_{\geq 0}$.  In order to show that
  $\mathfrak{c}$ and $\mathfrak{s}$ are inverses of each other we need
  to show that $\mathfrak
  s(\psi_{\underline{a}})=\phi_{\underline{a}}$, for
  all $n$-tuples $\underline{a}\in \mathbb{Z}^n_{\geq 0}$. 
  
  Let $\theta\in k[X_1,X_2, \cdots , X_d] \otimes k[t]$ be the
  polynomial $t-X_1^{a_1}X_2^{a_2}\dots X_n^{a_n}$. Then
  $p(\theta)= tI-x_1^{a_1}x_2^{a_2}\dots x_d^{a_d} \in \mathfrak
  g_1(A \otimes k[t])$, where $I$ is the identity matrix. Since
  $\det(x)=N_+(x)^2$ for any $x \in \mathfrak g_1(R)$, where $R$ is a
  $k$-algebra, we have
$$\det(p(\theta))=N_+(p(\theta))^2.$$ 
Using \eqref{section_AII} we also have
$$N_+(p(\theta))=\mathfrak s(\beta(\theta^{\otimes n})).$$
By the definition of the map $\beta$ we have
$\beta(\theta_a^{\otimes
  n})=\prod_{j=1}^n(t-\prod_{i=1}^db_j(y_i)^{a_i})$.

Hence, we get that
\begin{align*}
  \det(tI-x_1^{a_1}x_2^{a_2}\dots x_d^{a_d})
  =\mathfrak s(\prod_{j=1}^n(t-\prod_{i=1}^db_j(y_i)^{a_i}))^2.
 \end{align*}  
  Note that $\det(tI-x_1^{a_1}x_2^{a_2}\dots x_d^{a_d})$ is the
 characteristic polynomial of the matrix
 $x_1^{a_1}x_2^{a_2}\dots x_d^{a_d}$ and hence
 $\det(tI-x_1^{a_1}x_2^{a_2}\dots x_d^{a_d})=
 t^{2n}-Tr(x_1^{a_1}x_2^{a_2}\dots x_d^{a_d})t^{2n-1}+$ lower degree
 terms in $t$.
 On the other hand
 $$\mathfrak s(\prod_{j=1}^n(t-\prod_{i=1}^db_j(y_i)^{a_i}))^2
 =t^{2n}+2\sum_{j=1}^n\prod_{i}b_j(y_i)^{a_i} t^{2n+1}+\ \text{lower
   degree terms in}\ t.$$ Comparing the coefficients of $t^{2n-1}$ we
 get that $\mathfrak s(\mathfrak i(\phi_a))=\phi_a$.

  Hence, we get that the Chevalley restriction map
  $\mathfrak{i}:\mathfrak{c}^d//W_\mathfrak{c}\rightarrow
  \mathfrak{C}^d(\mathfrak{g}_1)//{G_0}$ is an isomorphism.

  \subsection{The symmetric pair \texorpdfstring{$AIII$}{}}
  Let $m,n$ be positive integers such that $m \geq n$. Let
  $\mathfrak{g}$ be the Lie algebra $\mathfrak{gl}_{n+m}(k)$ and let
  $\mathfrak{g}_0$ be the Lie algebra
 $$\left\{\begin{pmatrix}A&0\\ 0&B\end{pmatrix}:
   A \in \mathfrak{gl}_n(k), B\in \mathfrak{gl}_m(k)\right\},$$
 and let $\mathfrak{g}_1$ be the space 
 $$\left\{\begin{pmatrix}0&X\\ Y&0\end{pmatrix}:X, Y^t\in
   {\rm M}_{n\times m}(k)\right\}.$$ The pair
 $(\mathfrak{g}, \mathfrak{g}_0)$ is a symmetric pair for an
 involution with $\mathfrak{g}_0$, the space of invariants and
 $\mathfrak{g}_1$, the $-1$ eigenspace.  The adjoint action of the
 group $G_0={\rm GL}_n(k)\times {\rm GL}_m(k)$ on
 the space $\mathfrak{g}_1$ is given by
$$(g_1, g_2)(X, Y)=(g_1Xg_2^{-1}, g_2Yg_1^{-1}).$$
A Cartan subspace, denoted by $\mathfrak{c}$, of the
symmetric pair $(\mathfrak{g}, \mathfrak{g}_0)$ is given by the space
consisting of matrices of the form 
$$\begin{pmatrix}0&X\\X^t&0\end{pmatrix},$$ where $X=[\diag(b_1,
b_2,\dots, b_n), 0_{n\times m-n}]$ for some $b_i\in k$,
$i\in [n]$. The little Weyl group $W_\mathfrak{c}$ in this case is
isomorphic to $(\mathbb{Z}/2\mathbb{Z})^n\rtimes S_n$.

Let $\mathfrak{C}^d(\mathfrak{g}_1)$ be the commuting scheme
associated to the pair $(\mathfrak{g}, \mathfrak{g}_0)$. Let
$A=k[\mathfrak{C}^d(\mathfrak{g}_1)]$ be the coordinate ring of
$\mathfrak{C}^d(\mathfrak{g}_1)$ and let $(x_1,\dots, x_d)$ be the
universal point of $\mathfrak{C}^d(\mathfrak{g}_1)(A)$. We set
$$x_i=\begin{pmatrix}0&Q_i\\R_i&0\end{pmatrix}, Q_i, R_i^t\in
\mathfrak{gl}_{n\times m}(A).$$ Let $R$ be a $k$-algebra and let
$p:R[X_1,\dots, X_d]\rightarrow \mathfrak{gl}_{n+m}(A\otimes R)$ be
the map sending $X_i\mapsto x_i$ for all $i\in [d]$. Let
$R[X_1,\dots, X_d]^+$ be the algebra spanned by even degree
monomials. Note that the image of the map
$p: R[X_1,\dots, X_d]^+\rightarrow \mathfrak{gl}_{n+m}(A\otimes R)$ is
contained in
$\mathfrak{gl}_n(A\otimes R)\times \mathfrak{gl}_m(A\otimes R)$. Let
$q_1$ be the first projection of
$\mathfrak{gl}_n(A\otimes R)\times \mathfrak{gl}_m(A\otimes R)$.  Note
that the composite map
$$R[X_1,\dots, X_d]^+\xrightarrow{p} \mathfrak{gl}_n(A\otimes
R)\times \mathfrak{gl}_m(A\otimes
R)\xrightarrow{q_1}\mathfrak{gl}_n(A\otimes R) \xrightarrow{{\rm det}}
{A\otimes R}$$ has degree $n$. By Roby's theorem we get a map
$$\tilde{\mathfrak{s}}: TS^n(k[X_1,\dots, X_d]^+)\rightarrow A$$
such that
$\tilde{\mathfrak{s}}(\theta^{\otimes n})=\det\ \circ\ q_1\ \circ
p(\theta)$, for all $\theta\in k[X_1,\dots, X_d]^+$. Note that the
image of the map $\tilde{\mathfrak s}$ is contained in $A^{G_0}$.

The algebra $TS^nk[X_1,\dots, X_d]^+$ can be identified with the
subalgebra of $T^nk[X_1,\dots, X_d]$ of fixed points under
$W_{\mathfrak c}=S_n \ltimes (\mathbb Z/2\mathbb Z)^n$. Now, we give
an isomorphism of $TS^nk[X_1,\dots, X_d]^+$ with
$k[\mathfrak c^d]^{W_{\mathfrak c}}$. Let $(y_1,\dots, y_d)$ be the
tautological point of $\mathfrak{c}(k[\mathfrak c^d])$. Note that
$k[\mathfrak{c}^d]$ is a polynomial algebra in the variables
$b_j(y_i)$, $1 \leq i \leq d$ and $1 \leq j \leq n$. Let
$$\beta: T^nk[X_1,\dots, X_d] \rightarrow k[\mathfrak{c}^d]$$ be the
isomorphism of algebras defined by $\beta(X_{ij})=b_j(y_i)$, for
$i\in [d]$ and $j\in [n]$.  Here, $X_{ij}$ is the $i$-th variable in
the $j$-th copy of $k[X_1,\dots, X_d]$.  By restriction we have an
isomorphism of algebras
$$\beta: TS^nk[X_1,\dots, X_d]^+ \rightarrow k[\mathfrak c^d]^{W_{\mathfrak
c}}.$$ Then it follows that we have an algebra map
$$\mathfrak s:  k[\mathfrak c^d]^{W_{\mathfrak c}} \rightarrow
A^{G_0}$$ such that
$$\mathfrak{s}\beta(\theta^{\otimes n})=\det\circ q_1 \circ p(\theta),
\,\, \text{for all} \,\, \theta \in k[X_1,\dots, X_d]^+.$$ Thus we
obtain a map of schemes
$\mathfrak s: \mathfrak{C}^d(\mathfrak{g}_1)//G_0 \rightarrow
\mathfrak{c}^d//W_\mathfrak{c}$.
  
  Let $\mathfrak{i}:k[\mathfrak{C}^d(\mathfrak{g}_1)]^{G_0}\rightarrow
k[\mathfrak{c}^d]^{W_\mathfrak{c}}$ be the restriction map. We will show that
$\mathfrak s$ and $\mathfrak i$ are inverses of each other.

The ring of invariants $k[\mathfrak{C}^d(\mathfrak{g}_1)]^{G_0}$ is
generated by the images of the following polynomials, via the
restriction map
$k[\mathfrak{g}_1^d]^{G_0}\rightarrow
k[\mathfrak{C}^d(\mathfrak{g}_1)]^{G_0}$:
\begin{equation}\label{eq_AIII}
  {\rm Tr}(M_1M_2\dots M_k)\end{equation}
where $M_i=Q_{n_i}R_{m_i}$ for some $n_i, m_i\in [d]$; and we denote
the polynomial in \eqref{eq_AIII} by $P$ (see Section
\ref{invariants_gl_n_gl_m}).
We then have
$$\mathfrak{i}(P)=\sum_{j=1}^n\prod_{(n_i,
  m_i)}b_j(y_{n_i})b_j(y_{m_i}).$$
From the above equation, we get that
 $\mathfrak{i}(\phi_{\underline{a}})$ generate the $k$-algebra
 $k[\mathfrak{c}^d]^{W_\mathfrak{c}}$. Let $R$ be the algebra $k[t]$, and
consider the element $\theta=t-\prod_{(n_i, m_i)}X_{n_i}X_{m_i}$ in
$R[X_1, X_2,\dots, X_d]^+$. We then have
$$\beta(\theta^{\otimes n})=
\prod_{j=1}^n(t-\prod_{(n_i, m_i)}b_j(y_{n_i})b_j(y_{m_i})).$$
We then observe that $\det\circ q_1\circ p(\theta)$ is equal to
$$\det(t-\prod_{(n_i, m_i)}b_j(y_{n_i})b_j(y_{m_i})).$$
Since $\mathfrak{s}\circ \beta(\theta^{\otimes n})=\det\circ q_1\circ
p(\theta)$, we get that
$$t^n-\mathfrak{s}\left(\sum_{j=1}^n
  \prod_{(n_i, m_i)}b_j(y_{n_i})b_j(y_{m_i})\right)t^{n-1}+\cdots
=t^n-{\rm Tr}(\prod_{(n_i, m_i)}Q_{n_i}R_{m_i})t^{n-1}+\cdots$$
Comparing the coefficients of $t^{n-1}$ we get that $\mathfrak{s}\circ
\mathfrak{i}(P)=P$. Thus, we get that $\mathfrak{s}$ is a section
of the map $\mathfrak{i}$ and hence $\mathfrak i$ is an isomorphism.

\subsection{The symmetric pair \texorpdfstring{$BDI$}{}}
Let $m,n$ be positive integers such that $m \geq n$ and let
$\mathfrak{g}$ be the Lie algebra $\mathfrak{so}_{n+m}(k)$.  Let
$\mathfrak{g}_0$ be the Lie algebra
 $$\left\{\begin{pmatrix}A&0\\ 0&B\end{pmatrix}:A \in \mathfrak{so}_n(k), 
 B\in \mathfrak{so}_m(k)\right\},$$
 and let $\mathfrak{g}_1$ be the space 
 $$\left\{\begin{pmatrix}0&X\\ -X^t&0\end{pmatrix}:X\in
   {\rm M}_{n\times m}(k)\right\}.$$ The pair
 $(\mathfrak{g}, \mathfrak{g}_0)$ is a symmetric pair for an
 involution with $\mathfrak{g}_0$, the space of invariants and
 $\mathfrak{g}_1$, the $-1$ eigenspace.  The adjoint action of the
 group $G_0={\rm SO}_n(k)\times {\rm SO}_m(k)$ on
 the space $\mathfrak{g}_1$ is given by
$$(g_1, g_2)(X, Y)=(g_1Xg_2^{-1}, g_2Yg_1^{-1}).$$
 A Cartan subspace, denoted
by $\mathfrak{c}$, of the
symmetric pair $(\mathfrak{g}, \mathfrak{g}_0)$ is given by the space
consisting of matrices of the form 
$$\begin{pmatrix}0&X\\-X^t&0\end{pmatrix},$$ where 
$X$ is a matrix of the form $[A, 0_{n\times (m-n)}]$, with
$A=\diag(b_1, b_2,\dots, b_n)$ for some $b_i\in k$, $i\in [n]\}$.  
The little Weyl group $W_\mathfrak{c}$ is isomorphic to
$\mathbb{Z}/2\mathbb{Z}^n\rtimes S_n$. This case is similar to the
case of $AIII$, however since the space $\mathfrak{g}_1$ differs from
its counter part in $AIII$, we briefly describe the section for the
Chevalley restriction map.

Let $\mathfrak{C}^d(\mathfrak{g}_1)$ be the commuting scheme
associated to the pair $(\mathfrak{g}, \mathfrak{g}_0)$. Let
$A=k[\mathfrak{C}^d(\mathfrak{g}_1)]$, and let $(x_1,\dots, x_d)$ be
the universal point in $\mathfrak{C}^d(\mathfrak{g}_1)(A)$. We set
$$x_i=\begin{pmatrix}0&Q_i\\-Q_i^t&0\end{pmatrix}, Q_i\in
\mathfrak{gl}_{n\times m}(A).$$ Let $R$ be a $k$-algebra and let
$p:R[X_1,\dots, X_d]\rightarrow \mathfrak{gl}_{n+m}(A\otimes R)$ be
the map sending $X_i\mapsto x_i$ for all $i\in [d]$. Let
$R[X_1,\dots, X_d]^+$ be the subalgebra of $R[X_1,\dots, X_d]$ spanned
by the even degree monomials. Note that the image of the map
$p: R[X_1,\dots, X_d]^+\rightarrow \mathfrak{gl}_{n+m}(A\otimes R)$ is
contained in
$\mathfrak{gl}_n(A\otimes R)\times \mathfrak{gl}_m(A\otimes R)$. Let
$q_1$ be the first projection of
$\mathfrak{gl}_n(A\otimes R)\times \mathfrak{gl}_m(A\otimes R)$.  Note
that the composite map
$$R[X_1,\dots, X_d]^+\xrightarrow{p} \mathfrak{gl}_n(A\otimes
R)\times \mathfrak{gl}_m(A\otimes
R)\xrightarrow{q_1}\mathfrak{gl}_n(A\otimes R) \xrightarrow{{\rm det}}
{A\otimes R}$$ is multiplicative and has degree $n$. So by Roby's
theorem we get a map
$$\tilde{\mathfrak{s}}: TS^n(k[X_1,\dots, X_d]^+)\rightarrow A$$
such that
$\tilde{\mathfrak{s}}(\theta^{\otimes n})=\det\ \circ\ q_1\ \circ
p(\theta)$, for all $\theta\in k[X_1,\dots, X_d]^+$. We note that the
image of  $\tilde{\mathfrak{s}}$ belongs to $A^{G_0}$. 

Let $(y_1,\dots y_d)$ be a
tautological point of $k[\mathfrak{c}^d](B)$, where
$B=k[\mathfrak{c}^d]$. Let
$$\beta: TS^n(k[X_1,\dots, X_d]^+)\rightarrow k[\mathfrak{c}^d]$$
be the map $\beta(X_{ij})=b_j(y_i)$, for all $i\in [d]$ and
$j\in [n]$. Here, $X_{ij}$ be the variable $X_j$ in the $i$-th copy in
$TS^n(k[X_1,\dots, X_d])$. As in the case of $AIII$, $\beta$ restricts
to an isomorphism of algebras
$$\beta: TS^nk[X_1,\dots, X_d]^+ \rightarrow k[\mathfrak c^d]^{W_{\mathfrak
c}}.$$ Then it follows that we have an algebra map
$$\mathfrak s:  k[\mathfrak c^d]^{W_{\mathfrak c}} \rightarrow
A^{G_0}$$ such that
$$\mathfrak{s}\beta(\theta^{\otimes n})=\det\circ q_1 \circ p(\theta),
\,\, \text{for all} \,\, \theta \in k[X_1,\dots, X_d]^+.$$ Thus we
obtain a map of schemes
$\mathfrak s: \mathfrak{C}^d(\mathfrak{g}_1)//G_0 \rightarrow
\mathfrak{c}^d//W_\mathfrak{c}$.

We assume that $m$ is odd. The ring of invariants
$k[\mathfrak{C}^d(\mathfrak{g}_1)]^{G_0}$ is generated by the
images of the following polynomials, via the restriction map
$k[\mathfrak{g}_1^d]^{G_0}\rightarrow
k[\mathfrak{C}^d(\mathfrak{g}_1)]^{G_0}$:
\begin{equation}\label{gen_BDI}
  {\rm Tr}(M_1M_2\dots M_k)
  \end{equation}
  where $M_i=Q_{n_i}Q^t_{m_i}$ for some $n_i \in [d]$; and we
  denote the polynomial in \eqref{gen_BDI} by $P$ (see Subsection
  \ref{invariants_so_so}). We then have
$$\mathfrak{i}(P)=\sum_{j=1}^n\prod_{n_i}b_j(y_{n_i})b_i(y_{m_i}).$$
From the above equation, we get that
 $\mathfrak{i}(\phi_{\underline{a}})$ generate the $k$-algebra
 $k[\mathfrak{c}^d]^{W_\mathfrak{c}}$. 
Let $R$ be the algebra $k[t]$, and set $\theta=t-\prod_{n_i}
 X_{n_i}X_{m_i}$ in the algebra $R[X_1, X_2,\dots, X_d]^+$. We
then have
$$\beta(\theta^{\otimes n})=
\prod_{j=1}^n(t-\prod_{(n_i, m_i)}b_j(y_{n_i})b_j(y_{m_i})).$$
We then observe that $\det\circ q_1\circ p(\theta)$ is equal to
$$\det(t-\prod_{(n_i, m_i)}b_j(y_{n_i})b_j(y_{m_i})).$$
Since $\mathfrak{s}\circ \beta(\theta^{\otimes n})=\det\circ q_1\circ
p(\theta)$, we get that
$$t^n-\mathfrak{s}\left(\sum_{j=1}^n
  \prod_{(n_i, m_i)}b_j(y_{n_i})b_j(y_{m_i})\right)t^{n-1}+\cdots
=t^n-{\rm Tr}(\prod_{(n_i, m_i)}Q_{n_i}Q^t_{m_i})t^{n-1}+\cdots$$
Comparing the coefficients of $t^{n-1}$ we get that $\mathfrak{s}\circ
\mathfrak{\mathfrak{i}}(P)=P$. Thus, we get that $\mathfrak{s}$ is a section
of the map $\mathfrak{i}$.

\begin{remark}\normalfont
When $m=n$ is even, the ring of invariants
$k[\mathfrak{C}_d(\mathfrak{g}_1)]^{G_0}$ has additional invariants of the form
\begin{equation}\label{new_gen_BDI}
 \det(T_1\otimes A_1+T_2\otimes A_2+\cdots+T_d\otimes A_d)
\end{equation}
where $T=(T_1,\dots, T_d)$ is an element of $M^d_{r\times r}$ and 
$r\geq 1$. We denote by $\psi_T$ the invariant in
\eqref{new_gen_BDI}. Consider the element
$$\theta=\det(T_1\otimes X_1I_{n\times n}+T_2\otimes X_2I_{n\times
  n}+\cdots+T_d\otimes X_dI_{n\times n})$$
in the ring $R[X_1,\dots, X_d]^+$, where $R=k[M^d_{r\times r}]$. The
identity $\mathfrak{s}\circ \beta(\theta^{\otimes n})=\det\circ
q_1\circ p(\theta)$ implies that $\mathfrak{s}\circ
\mathfrak{i}(\psi_T)=\psi_T$. Now, under the assumption that the
invariants in \eqref{new_gen_BDI} and \eqref{gen_BDI} together generate the
algebra $k[\mathfrak{C}_d(\mathfrak{g}_1)]^{G_0}$, we get that the map
$\mathfrak{s}$ is the section of the Chevalley restriction map. 
\end{remark}

\subsection{The symmetric pair \texorpdfstring{$CI$}{}}  
  Let $V$ be a $2n$-dimensional $k$-vector space and let $\omega$
  be a non-degenerate skew symmetric bilinear form on $V$. Let $(v_1,
  v_2,\dots v_n, v_{-1}, \dots, v_{-n})$ be a Witt basis for the pair
  $(V, \omega)$. The matrix of the form $\omega$ is equal to
  $$J=\begin{pmatrix}0&I_n\\I_n&0\end{pmatrix},$$
  where $I_n$ is the $n\times n$ identity matrix. Let $V_+$ and
  $V_{-}$ be the spaces $\langle v_1,\dots, v_n\rangle$ and
  $\langle v_{-1},\dots, v_{-n}\rangle$ respectively.  In this basis
  we identify $\mathfrak{g}=\mathfrak{sp}(V,\omega)$ as a Lie
  subalgebra of $\mathfrak{gl}_{2n}(k)$. Let $T \in GL(V)$ be such that
  $T|_{V_+}=id$ and $T|_{V_-}=-id$. Then the conjugation by $T$
  defines an involution of $\mathfrak g$ such that the Lie algebra
  $\mathfrak g_0$ is identified with the Lie
  algebra
  $$\{\begin{pmatrix}A&0\\0&-A^t\end{pmatrix}:A\in \mathfrak{gl}_n(k)\},$$
  and $\mathfrak{g}_1$ is the space 
  $$\{\begin{pmatrix}0&X\\Y&0\end{pmatrix}:X, Y\in \mathfrak{gl}_n(k),
  X^t=X, Y^t=Y\}.$$ Let $G_0$ be the connected subgroup of
  ${\rm Sp}(V, \omega)$ with $\mathfrak{g}_0$ as its Lie algebra. We
  identify $G_0$ with ${\rm GL}_n(k)$. As a $G_0$ module $\mathfrak{g}_1$ is
  isomorphic to ${\rm Sym}^2(k^n) \oplus ({\rm Sym}^2(k^n))^*$.
  We choose the Cartan subspace $\mathfrak{c}$ to be the subspace
  $$\{\begin{pmatrix}0&X\\X&0\end{pmatrix}:X=\diag(b_1,\dots, b_n),
b_i\in k, i\in [n]\}$$ of $\mathfrak{g}_1$. Then the little Weyl
group, denoted by $W_{\mathfrak c}$, is
isomorphic to $S_n \ltimes (\mathbb Z/2\mathbb Z)^n$.
 
  Let $\mathfrak{C}^d(\mathfrak{g}_1)$ be the $d$-fold commuting
  scheme associated with the symmetric pair
  $(\mathfrak g, \mathfrak g_0)$ and let
  $A=k[\mathfrak{C}^d(\mathfrak{g}_1)]$ be the coordinate ring. Let
  $(x_1,\dots, x_d)$ be the universal point in
  $\mathfrak{C}^d(\mathfrak{g}_1)(A)$. Let $Q_i$ and $R_i$ be matrices
  with values in $A$ such that
  \begin{equation}\label{uni_CI}
    x_i=\begin{pmatrix}0&Q_i\\R_i&0\end{pmatrix}, i \in [d].
    \end{equation}
  Note that $Q_i^t=Q_i$ and $R_i^t=R_i$, for all $i\in [d]$. 
  Let $R$ be a $k$ algebra and let $p:R[X_1,\dots, X_d]\rightarrow
  \mathfrak{gl}_{2n}(A \otimes R)$
  be the map defined by $X_i\mapsto x_i$. Let $R[X_1,\dots, X_d]^+$
  be the $R$-subalgebra generated by even degree monomials. Note that
  the image of the map $p$ is contained in $\mathfrak{gl}_n (A
  \otimes R) \times \mathfrak{gl}_n(A \otimes R)$. Let $q_1$ and $q_2$
  be the first and second projections
  of $\mathfrak{gl}_n (A
  \otimes R) \times \mathfrak{gl}_n(A \otimes R)$ respectively. We have
  $$R[X_1,\dots, X_d]^+\xrightarrow{p} \mathfrak{gl}_n (A
  \otimes R) \times \mathfrak{gl}_n(A \otimes R)\xrightarrow{det\
    \circ\ q_1}A\otimes R.$$
  
  Since the map $\det$ is multiplicative, the composition $\det\circ q_1 \circ p$ is
multiplicative and homogeneous of degree $n$. Thus by Roby's theorem
we get a homomorphism of $k$-algebras:
$$\tilde{\mathfrak s}:TS^n(k[X_1,\dots,X_d]^+)\rightarrow A$$
such that 
\begin{equation}
  \tilde{\mathfrak s}(\theta^{\otimes n})=
det \circ q_1(p(\theta)), \,\, \text{for all} \,\, \theta \in k[X_1,\dots,
X_d].
\end{equation}
Note that the image of the map $\tilde{\mathfrak s}$ is contained in
$A^{G_0}$.

The subalgebra $TS^nk[X_1,\dots, X_d]^+$ can be identified with the
subalgebra of $T^nk[X_1,\dots, X_d]$ of fixed points under
$S_n \ltimes (\mathbb Z/2\mathbb Z)^n$. Now, we give an isomorphism of
$TS^nk[X_1,\dots, X_d]^+$ with $k[\mathfrak c^d]^{W_{\mathfrak
    c}}$. Let $(y_1,\dots, y_d)$ be the tautological point of
$\mathfrak{c}(k[\mathfrak c^d])$. Note that $k[\mathfrak{c}^d]$ is a
polynomial algebra in the variables $b_j(y_i)$, $1 \leq i \leq d$ and
$1 \leq j \leq n$. Let
$$\beta: T^nk[X_1,\dots, X_d] \rightarrow k[\mathfrak{c}^d]$$ be the
isomorphism of algebras defined by $\beta(X_{ij})=b_j(y_i)$, for
$i\in [d]$ and $j\in [n]$.  Here, $X_{ij}$ is the $i$-th variable in
the $j$-th copy of $k[X_1,\dots, X_d]$.  By restriction we have an
isomorphism of algebras
$$\beta: TS^nk[X_1,\dots, X_d]^+ \rightarrow k[\mathfrak c^d]^{W_{\mathfrak
c}}.$$ Then it follows that we have an algebra map
$$\mathfrak s:  k[\mathfrak c^d]^{W_{\mathfrak c}} \rightarrow
A^{G_0}$$ such that
$$\mathfrak{s}\beta(\theta^{\otimes n})=\det\circ q_1 \circ p(\theta),
\,\, \text{for all} \,\, \theta \in k[X_1,\dots, X_d]^+.$$ Thus we
obtain a map of schemes
$\mathfrak s: \mathfrak{C}^d(\mathfrak{g}_1)//G_0 \rightarrow
\mathfrak{c}^d//W_\mathfrak{c}$.
  
  Let $\mathfrak{i}:k[\mathfrak{C}^d(\mathfrak{g}_1)]^{G_0}\rightarrow
k[\mathfrak{c}^d]^{W_\mathfrak{c}}$ be the Chevalley restriction map. We will show that
$\mathfrak s$ and $\mathfrak i$ are inverses of each other.
  
Recall that $\mathfrak{g}_1$ as a $G_0$-module is isomorphic to
${\rm Sym}^2(k^n) \oplus ({\rm Sym}^2(k^n))^*$.  The ring of
invariants $k[\mathfrak{C}^d(\mathfrak{g}_1)]^{G_0}$ is generated
by the images of the following polynomials, via the restriction map
$k[\mathfrak{g}_1^d]^{G_0}\rightarrow
k[\mathfrak{C}^d(\mathfrak{g}_1)]^{G_0}$:
\begin{equation}\label{gen_AIII}
  {\rm Tr}(M_1M_2\dots M_k)\end{equation}
where $M_i=Q_{n_i}R_{m_i}$ for some $n_i, m_i\in [d]$; and we denote
the polynomial in \eqref{gen_AIII} by $P$ (see Subsection
\ref{invariants_sym_gl_n}). Here $Q_{n_i}$ and $R_{n_i}$ are as
defined in  \eqref{uni_CI}. 
We then have
$$\mathfrak{i}(P)=\sum_{j=1}^n\prod_{(n_i,
  m_i)}b_j(y_{n_i})b_j(y_{m_i}).$$ From the above equation, we get
that $\mathfrak{i}(\phi_{\underline{a}})$ generate the $k$-algebra
$k[\mathfrak{c}^d]^{W_\mathfrak{c}}$.  Let $R$ be the algebra $k[t]$,
and consider the element $\theta=t-\prod_{(n_i, m_i)}X_{n_i}X_{m_i}$
in $R[X_1, X_2,\dots, X_d]^+$. We then have
$$\beta(\theta^{\otimes n})=
\prod_{j=1}^n(t-\prod_{(n_i, m_i)}b_j(y_{n_i})b_j(y_{m_i})).$$
We then observe that $\det\circ q_1\circ p(\theta)$ is equal to
$$\det(t-\prod_{(n_i, m_i)}Q_{n_i}R_{m_i}).$$
Since $\mathfrak{s}\circ \beta(\theta^{\otimes n})=\det\circ q_1\circ
p(\theta)$, we get that
$$t^n-\mathfrak{s}\left(\sum_{j=1}^n
  \prod_{(n_i, m_i)}b_j(y_{n_i})b_j(y_{m_i})\right)t^{n-1}+\cdots
=t^n-{\rm Tr}(\prod_{(n_i, m_i)}Q_{n_i}R_{m_i})t^{n-1}+\cdots$$
Comparing the coefficients of $t^{n-1}$ we get that $\mathfrak{s}\circ
\mathfrak{i}(P)=P$. Thus, we get that $\mathfrak{s}$ is a section
of the map $\mathfrak{i}$ and hence $\mathfrak i$ is an isomorphism.

\section{The symmetric pair of type \texorpdfstring{$CII$}{}}
Let $n$ be a positive integer and let $n=q+r$ for some positive
integers $q$ and $r$ with $r\leq q$. Let $V$ be a $2n$ dimensional
$k$-vector space and let $\omega$ be a non-degenerate symplectic form
on $V$. Let
$$(w_1, w_{-1}, w_{2}, w_{-2}, \dots, w_{n}, w_{-n})$$
be a Witt-basis for $V$ such that $\omega(w_i, w_j)=1$, for $i+j=0$,
$i>0$ and $\omega(w_i, w_j)=0$, for all $i, j$ such that $i+j\neq
0$. Let $V_1$ and $V_2$ be the subspaces of $V$ spanned by
$\{w_{\pm 1}, w_{\pm 2},\dots, w_{\pm r}\}$ and
$\{w_{\pm (r+1)}, w_{\pm (r+2)},\dots, w_{\pm (q+r)}\}$
respectively. Let $T_s$ be the $2s\times 2s$ matrix given by
$\diag(\mu, \mu,\dots, \mu)$, where $\mu$ is the
matrix $$\begin{pmatrix}0&1\\-1&0\end{pmatrix}.$$
Let $\mathfrak{g}$ be the Lie algebra defined by
$$\{X\in \End_k(V):\omega(Xv,w)+\omega(v, Xw)=0\}.$$
Using the above Witt-basis, we identify $\mathfrak{g}$ with the Lie
algebra $\mathfrak{sp}_{2n}(k)$. Let $\mathfrak{g}_0$ be the Lie
subalgebra of $\mathfrak{g}$ consisting of matrices
$$\begin{pmatrix}A&0\\0&B\end{pmatrix}$$
where $A$ and $B$ belong to $\mathfrak{sp}_{2r}(k)$ and
$\mathfrak{sp}_{2q}(k)$ respectively. Let $\mathfrak{g}_1$ be the 
subspace of $\mathfrak{g}$ consisting of matrices of the form
$$\begin{pmatrix}
  0&X\\T_rX^tT_q&0
\end{pmatrix}.$$ Note that $(\mathfrak{g}, \mathfrak{g}_0)$ is a
symmetric pair of type $CII$.  The Cartan subspace $\mathfrak{c}$ of
$(\mathfrak{g}, \mathfrak{g}_0)$ is given by the set of matrices of
the form 
$$\begin{pmatrix}
  0&X\\T_rX^tT_q&0
\end{pmatrix}.$$ where $X$ is of the form
$\begin{pmatrix}BT_r\\0_{q-r}\end{pmatrix}$ and $B$ is the diagonal
matrix $[b_1, b_1,\dots, b_r, b_r]$. The little Weyl group
$W_\mathfrak{c}$ is equal to $(\mathbb{Z}/2\mathbb{Z})^r\rtimes S_r$.

Let $\mathfrak{C}^d(\mathfrak{g}_1)$ be the $d$-fold commuting scheme
attached with the pair $(\mathfrak{g}, \mathfrak{g}_0)$. Let $A$ be
the coordinate ring of the affine scheme
$\mathfrak{C}^d(\mathfrak{g}_1)$, and let $(x_1, x_2,\dots, x_d)$ be
the universal point of the scheme $\mathfrak{C}^d(\mathfrak{g}_1)$. We
then have
$$x_i=\begin{pmatrix}0& X_i\\ T_rX_i^tT_q&0\end{pmatrix},$$
for $i\in [d]$. Let $R$ be a $k$-algebra, and let
$p:R[X_1,\dots, X_d]\rightarrow \mathfrak{gl}_{2n}(A\otimes R)$ be the
map $X_i\mapsto x_i$, for some $i\in [d]$. Let
$R[X_1, X_2,\dots, X_d]^+$ be the subalgebra consisting of
even degree polynomials in $R[X_1, X_2,\dots, X_d]$. We note that the
image of $p$ restricted to $R[X_1, X_2,\dots, X_d]^+$ lands in
$\mathfrak{gl}_{2r}(A\otimes R)\times \mathfrak{gl}_{2q}(A\otimes
R)$. Moreover, from the relations
$$X_iT_rX_j^t=X_jT_rX_i^t$$
and
$$X_i^tT_qX_j=X_j^tT_qX_i$$
we get that the image of the composite map
$$R[X_1,\dots, X_d]^+\xrightarrow{p}\mathfrak{gl}_{2r}(A\otimes
R)\times \mathfrak{gl}_{2q}(A\otimes R)$$
is contained in $\mathfrak{g}^+_{r}(A\otimes R)\times
\mathfrak{g}^+_q(A\otimes R)$. Here, $\mathfrak{g}_m^+(A\otimes R)$ is
given by 
$$\{X\in M_{2m\times 2m}(A\otimes R):T_mX-X^tT_m=0\}.$$
Consider the following composite map 
$$R[X_1, X_2,\dots, X_d]^+\xrightarrow{p} \mathfrak{g}^{+}_{r}(A\otimes
R)\times \mathfrak{g}^{+}_q(A\otimes R)\xrightarrow{q_1}
\mathfrak{g}^+_{r}(A\otimes R)\xrightarrow{N_+}A\otimes R$$ Here,
$N_+$ is the Pfaffian norm map as defined in \ref{paffian_map_AII} and
$q_1$ is the first projection.  Note that the map
$N_+\circ q_1\circ p$ is a degree $r$ map and hence Roby's theorem
implies that there exists an $k$-algebra homomorphism
$$\tilde{s}:TS^rR[X_1, \dots, X_d]^+\rightarrow A^{G_0}$$
such that $\tilde{s}(\theta^{\otimes r})= N_+\circ q_1\circ p(\theta)$
for all $\theta\in R[X_1,\dots, X_d]^+$.

Let $(y_1,\dots, y_d)$ be the tautological point of
$\mathfrak{c}(k[\mathfrak c^d])$. Note that $k[\mathfrak{c}^d]$ is a
polynomial algebra in the variables $b_j(y_i)$, $1 \leq i \leq d$ and
$1 \leq j \leq r$. Let
$$\beta: T^rk[X_1,\dots, X_d] \rightarrow k[\mathfrak{c}^d]$$ be the
isomorphism of algebras defined by $\beta(X_{ij})=b_j(y_i)$, for
$i\in [d]$ and $j\in [r]$.  Here, $X_{ij}$ is the $i$-th variable in
the $j$-th copy of $k[X_1,\dots, X_d]$.  By restriction we have an
isomorphism of algebras
$$\beta: TS^rk[X_1,\dots, X_d]^+ \rightarrow k[\mathfrak c^d]^{W_{\mathfrak
c}}.$$ Then it follows that we have an algebra map
$$\mathfrak s:  k[\mathfrak c^d]^{W_{\mathfrak c}} \rightarrow
A^{G_0}$$ such that
$$\mathfrak{s}\beta(\theta^{\otimes r})=N_+\circ q_1 \circ p(\theta),
\,\, \text{for all} \,\, \theta \in k[X_1,\dots, X_d]^+.$$

Note that $\mathfrak{g}_1$ as a $G_0$ module is isomorphic to
$M_{2r\times 2q}(k)$. Since $\mathfrak{C}^d(\mathfrak{g}_1)$ is a closed 
subscheme of $\mathfrak g_1^d$, the ring of invariants
$k[\mathfrak{C}^d(\mathfrak{g}_1)]^{G_0}$ is generated by polynomials
of the form
$${\rm Tr}(M_{n_1m_1}M_{n_2m_2}\dots M_{n_lm_l})$$
where $M_{ij}=M_i^tT_rM_jT_q$ and $n_i, m_i\in [d]$ and
$M_i\in M_{2r\times 2q}$ is a $2r\times 2q$ matrix (see subsection
\ref{invariants-type-cii}). We denote the above polynomial by $P$. Let
$\mathfrak{i}: k[\mathfrak{C}^d(\mathfrak{g}_1)]^{G_0} \rightarrow
k[\mathfrak{c}^d]^{W_\mathfrak{c}}$ be the Chevalley's restriction
map. Note that $\mathfrak{i}(P)$ is given by
 $$\mathfrak{i}(P)=2\sum_{j=1}^p\prod_{(n_i,
   m_i)}b_j(y_{n_i})b_{j}(y_{m_i}).$$
 Hence, by the explicit nature of the polynomials $\mathfrak{i}(P)$,
 we get that $k[\mathfrak{c}^d]^{W_\mathfrak{c}}$ is generated by
 polynomials of the form $\mathfrak{i}(P)$. Let $R$ be the algebra
 $k[t]$, and consider the element $\theta=t-\prod_{(n_i,
   m_i)}X_{n_i}X_{m_i}$ in $R[X_1, X_2,\dots, X_d]^+$. We then have
 \begin{equation}\label{CII_beta}
   \beta(\theta^{\otimes r})=\prod_{j=1}
   \left(t-\prod_{(n_i, m_i)}b_j(y_{n_i})b_j(y_{m_i})\right)
   \end{equation}
 Note that $N_+\circ q_1\circ p(\theta)$ is equal to
 \begin{equation}\label{CII_norm}
   N_+\left(t-\prod_{(n_i, m_i)}X_{n_i}T_rX^t_{m_i}T_q\right)
 \end{equation}
 Taking square of \eqref{CII_beta} and \eqref{CII_norm}  and using
 $\mathfrak{s}\circ \beta(\theta^{\otimes r})=N_+\circ q_1\circ
 p(\theta)$ we get that
 $$\det\left(t-\prod_{(n_i, m_i)}X_{n_i}T_rX^t_{m_i}T_q\right)
 =\prod_{j=1} \left(t-\prod_{(n_i,
     m_i)}b_j(y_{n_i})b_j(y_{m_i})\right)^2.$$ Comparing the
 coefficients of $t^{2r-1}$, we get that
 $\mathfrak{s}\circ \mathfrak{i}(P)=P$, for all
 $P\in k[\mathfrak{C}^d(\mathfrak{g}_1)]^{G_0}$. Hence $\mathfrak{s}$
 is a section of $\mathfrak{i}$. 

 This concludes the proof of Theorem \ref{main}.
  \begin{remark}\normalfont\label{except}
    The above techniques do not work in the case where
    ($\mathfrak{g}, \mathfrak{g}_0)$ is equal to
    $(\mathfrak{so}_{2n}(k), \mathfrak{gl}_n(k))$. Note that the rank
    of the above symmetric pair is $[n/2]$ and the little Weyl group
    is isomorphic to the Weyl group of type $B_n$.  Let
    $(x_1, x_2,\dots, x_d)$ be the universal point of
    $\mathfrak{C}_d(\mathfrak{g}_1)(A)$, where
    $A=k[\mathfrak{C}_d(\mathfrak{g}_1)]$. For any $k$-algebra $R$,
    let
    $p:R[X_1, X_2,\dots, X_d]^+\rightarrow \mathfrak{gl}_{2n}(A\otimes
    R)$ be the map $X_i\mapsto x_i$. The image of the above map lands
    in $\mathfrak{g}_0(A\otimes R)=\mathfrak{gl}_{n}(A\otimes
    R)$. Since the determinant has degree $n$, we do not know any
    multiplicative map on the image of $p$ of degree $[n/2]$.
    \end{remark}

    {\bf Acknowledgements:} We thank Prof. Michael Brion for carefully
    reading the manuscipt and for his comments which helped us to
    improve the presentation of the article.

\bibliography{./biblio} \bibliographystyle{amsalpha}
Department of Mathematics and Statistics,\\
    Indian Institute of Technology Kanpur,\\
    U.P. India, 208016.\\
    email:\ \texttt{nsantosh@iitk.ac.in}; \texttt{santosha@iitk.ac.in}
  \end{document}